\documentclass[11pt]{article}

\usepackage{amsmath}
\usepackage{amsthm}
\usepackage{amsfonts}
\usepackage{eucal}
\usepackage[dvips]{epsfig}
\psdraft
\psfull
\usepackage{calrsfs}
\usepackage{amssymb}

\def\CC{{\mathcal C}}

\def\Z{{\bf Z}}

\def\F{{\bf F}}

\def\N{{\bf N}}

\def\OO{{\mathcal O}}

\def\car{\mathop{\rm car}\nolimits}
\def\ord{\mathop{\rm ord}\nolimits}

\def\Aut{\mathop{\rm Aut}\nolimits}

\def\PGL{\mathop{\rm PGL}\nolimits}

\def\blok{\hskip .5em \vrule width 1.0ex height 1.0ex depth 0.3ex}

\newtheorem{theorem}{Th\'eor\`eme}[section]
\newtheorem{theo}[theorem]{Th\'eor\`eme}
\newtheorem{lem}[theorem]{Lemme}
\newtheorem{cor}[theorem]{Corollaire}
\newtheorem{prop}[theorem]{Proposition}

\theoremstyle{definition}

\newtheorem{example}[theorem]{Exemple}

\theoremstyle{remark}
\newtheorem{rem}[theorem]{Remarque}

\newif\ifnormalesBeweisEnde

  {\vskip 0.3ex plus 0.5ex minus 0ex \pagebreak[1]
   \global\normalesBeweisEndetrue
   \trivlist
   \item[\hskip\labelsep {\textsc{Proof} \rm of #1}:]}%
  {\ifnormalesBeweisEnde \EndOfBeweis \fi
   \endtrivlist
   \vskip 1ex plus 1ex minus 0ex \pagebreak[2]}
\newenvironment{pr}
  {\vskip 0.3ex plus 0.5ex minus 0ex \pagebreak[1]
   \global\normalesBeweisEndetrue
   \trivlist
   \item[\hskip\labelsep \textsc{Preuve} :]}%
  {\ifnormalesBeweisEnde \EndOfBeweis \fi
   \endtrivlist
   \vskip 1ex plus 1ex minus 0ex \pagebreak[2]}
  {\ifnormalesBeweisEnde \EndOfBeweis \fi
   \endtrivlist
   \vskip 1ex plus 1ex minus 0ex \pagebreak[2]}
\def\EndOfBeweis{\hskip .5em \vrule width 1.0ex height 1.0ex depth 0.3ex}

\begin{document}

\title{Rel\`evements des rev\^etements de courbes \\ faiblement ramifi\'es}

\author{Gunther Cornelissen \\ Universiteit Utrecht, Mathematisch Instituut \\  Postbus 80.010, 3508 TA Utrecht, Nederland \\ cornelis@math.uu.nl \\ \\ et \\  \\ Ariane M\'ezard \\ Universit\'e Paris-Sud \\ 91405 Orsay Cedex, France \\ ariane.mezard@math.u-psud.fr}

\date{}
\maketitle

\begin{quote}

\noindent
{\bf Abstract. } Let $X$ be a smooth projective curve over a perfect field of
characteristic  $p>0$ and $G$ a finite group of automorphism of $X$. Let
$\nu(X,G)$ be the characteristic of the versal 
equivariant deformation
ring $R(X,G)$ of $(X,G)$. When the ramification is weak
(i.e., all second ramification groups are trivial),
 we prove that $\nu(X,G)\in\{0,p\}$ and we compute
$R(X,G)$.\\

\noindent
{\bf R\'esum\'e.} Soit $X$ une courbe projective lisse sur un corps parfait de
caract\'eristique $p>0$ et $G$ un groupe fini d'automorphismes de $X$. Nous
consid\'erons la caract\'eristique $\nu(X,G)$ de l'anneau versel $R(X,G)$
de d\'eformations
\'equivariantes de $(X,G)$. Dans le cas d'une ramification faible
(o\`u tous les seconds groupes de ramification sont triviaux), nous d\'emontrons
que $\nu(X,G)\in\{0,p\}$ et nous calculons $R(X,G)$.\\
\end{quote}

\section{Introduction}
Soit $X$ une courbe alg\'ebrique projective lisse de genre $g$ sur un corps
parfait $k$ de caract\'eristique $p>0$ munie de l'action d'un groupe fini
d'automorphismes $G\subset\Aut_k X$. Soit $W(k)$ l'anneau des vecteurs de Witt
de $k$ et $\widehat{\CC}_k$ la cat\'egorie des $W(k)$-alg\`ebres locales
noeth\'eriennes compl\`etes de corps r\'esiduel $k$. Soit $R$ un objet de $\widehat{\CC}_k$.
Un rel\`evement \'equivariant du couple $(X,G)$ \`a $R$ est une courbe 
$\tilde{X}$ propre et lisse sur $R$ telle que $G\subset\Aut_R \tilde{X}$
et telle que l'isomorphisme $\tilde{X}\times_Rk\simeq X$ soit
$G$-\'equivariant. Si la ramification est mod\'er\'ee, il
existe un rel\`evement de $(X,G)$ \`a un anneau de $\widehat{\CC}_k$ de
caract\'eristique z\'ero (\cite{Gr}). Lorsque la ramification est sauvage, 
la situation est beaucoup plus myst\'erieuse : il
n'existe pas forc\'ement de rel\`evement de $(X,G)$ en caract\'eristique 
z\'ero.\\

\noindent
{\bf Exemple.} Supposons $p\geq 3.$ Consid\'erons la courbe $X$ sur
$k=\overline{\F}_p$ normalis\'ee de la courbe d'\'equation $(x^p-x)(y^p-y)=1$.
La courbe $X$ a pour genre $g=(p-1)^2$.
Soit $G=(\Z/p\Z)^2\rtimes D_{p-1}$ o\`u $D_n$ d\'esigne le groupe di\'edral
d'ordre $2n$. Nous avons $G\subset\Aut_k X$. Un rel\`evement $(\tilde{X},G)$ de $(X,G)$ en
caract\'eristique z\'ero aurait donc au moins $|G|$ automorphismes. Ceci est
exclu pour $p\geq 41$ car, d'apr\`es la borne de Hurwitz, $|\Aut \tilde{X}|\leq
84(g-1)$.\\

\noindent
Une question naturelle se pose alors : le couple $(X,G)$ se rel\`eve-t-il
\`a un objet de $\widehat{\CC}_k$ de caract\'eristique $p^2$ ?  Plus 
g\'en\'eralement, quelle est la plus grande puissance $p^n$ de $p$ pour
laquelle il existe un rel\`evement de $(X,G)$ \`a un objet $R$ de
$\widehat{\CC}_k$
de caract\'eristique $p^n$ ?
Traduisons ces questions dans le formalisme de la th\'eorie des
d\'eformations : soit $R(X,G)$ l'anneau versel de 
d\'eformations \'equivariantes
de $(X,G)$. Si $W(k)\subset R(X,G)$, 
nous posons $\nu(X,G)=0$. Sinon nous posons 
$\nu(X,G)=p^{n+1}$ 
o\`u $n\in\N$ est le plus grand entier
tel que $p^n\neq 0$ dans
$R(X,G)$. Il s'agit donc de calculer $\nu(X,G)$.

\begin{rem}
L'entier $\nu(X,G)$ est
la caract\'eristique de $R(X,G)$ au sens suivant : un anneau $A$ non nul 
d'unit\'e $1$ a pour caract\'eristique $n 
\in \N$ si $n$ est un g\'en\'erateur positif du noyau de l'unique homomorphisme 
 $\Z \rightarrow A: 1 \mapsto 1$ (voir \cite{Bourbaki} \S 8, no.\ 8). \\
\end{rem}

\noindent L'objet de cet article est le calcul de $\nu(X,G)$ pour une action faiblement
ramifi\'ee au sens suivant : soit $P$ un point de ramification de la courbe
$X$.
Soit $G_P$ le stabilisateur de $P$ d'uniformisante $x$. Les groupes de
ramification de $G_P$ sont not\'es
$$G_{P,0}=G_P=\{\sigma\in G,\sigma P=P\},\;\;\; G_{P,i}=\{\sigma\in
G,\ord_p(\sigma x-x)>i\}, \;\; i>0.$$
L'action de $G$ sur $X$ est dite {\it faiblement ramifi\'ee} si tous les
groupes de ramification $G_{P,2}$ sont triviaux. Notons que la ramification
mod\'er\'ee (i.e.\ $G_{P,1}=0$) implique la ramification faible. Une action
faible est donc une action sauvage la ``moins ramifi\'ee possible''.\\
Si l'action de $G$ sur $X$ est faiblement ramifi\'ee alors pour tout point $P$ de ramification de la courbe,
il existe $t\geq 0$ et $n$ divisant $p^t-1$ tels que
$$G_P=(\Z/p\Z)^t\rtimes\Z/n\Z,$$
o\`u $(\Z/p\Z)^t$ est identifi\'e \`a $\F_q$ pour $q=p^t$, l'action semi-directe est donn\'ee par $\mu_n$ sur ${\bf G}_a$ (voir
\cite{CoKa2} \S 1, \cite{Se} VI.2) et $G_P$ est de conducteur $1$. 
Remarquons que l'action d'un groupe fini sur une courbe ordinaire est toujours
faiblement ramifi\'ee (voir \cite{Na}).
Le r\'esultat principal de cet article est le
th\'eor\`eme suivant :
\begin{theo}
\label{theo1}
Soit $X$ une courbe projective lisse d\'efinie sur un corps parfait de
caract\'eristique $p>0$ et munie d'une action faiblement ramifi\'ee de groupe
$G$. Alors l'une des propri\'et\'es suivantes est satisfaite\;:\\ 
\textup{a.} $R(X,G)$ est annul\'e par $p$ ;\\
\textup{b.} les groupes de ramification de $G$ non triviaux
d'ordre divisible par $p$ sont de l'une des trois formes suivantes :\\
\indent \textup{i.} groupe cyclique $\Z/p\Z$ d'ordre $p$ ;\\
\indent \textup{ii.} groupe di\'edral $D_p=\Z/p\Z\rtimes\Z/2\Z$ d'ordre $2p$ si $p>2$ \\
\indent ou $(\Z/2\Z)^2$ si $p=2$ ;\\
\indent \textup{iii.} groupe t\'etra\'edrique $A_4=(\Z/2\Z)^2\rtimes\Z/3\Z$ pour $p=2$.\\
Dans le cas \textup{b.}, l'anneau $R(X,G)$ est d'intersection compl\`ete
relative sur $W(k)$ ; 
en particulier $R(X,G)$ est plat sur $W(k)$.
\end{theo}

\noindent
Nous obtenons ainsi $\nu(X,G) \in \{0,p\}$ pour une action faiblement ramifi\'ee.\\

\noindent Le plan de cet article est le suivant. La seconde partie traduit le r\'esultat principal en un th\'eor\`eme ``local'' (\ref{thlocal}). La troisi\`eme partie contient la d\'emon\-stra\-tion de ce th\'eor\`eme ``local''.  Nous avons report\'e dans la quatri\`eme partie les cons\'equences du th\'eor\`eme \ref{theo1} et des questions ouvertes
qu'il suscite. En particulier, en corollaire du th\'eor\`eme
\ref{theo1}, nous menons \`a bien le calcul des anneaux universels de d\'eformations lorsque l'action est faiblement
ramifi\'ee (corollaire \ref{coroloc}). \\

{\footnotesize La collaboration entre les deux auteurs a \'et\'e initi\'ee lors d'un s\'ejour
tr\`es enrichissant de A. M\'ezard \`a l'universit\'e de Leiden. A. M\'ezard
remercie B. de Smit et B. Edixhoven de leur accueil. Les auteurs remercient le rapporteur dont les remarques ont permis 
de clarifier la r\'edaction de cet article.}\\
\section{Du global au local}
\label{locglo}
Soit $k$ un corps parfait de caract\'eristique $p>0$.
Soit $\CC_k$ la cat\'egorie des $W(k)$-alg\`ebres locales artiniennes de corps
r\'esiduel $k$. Ainsi $\CC_k$ est une sous-cat\'egorie pleine de 
la cat\'egorie $\widehat{\CC}_k$ des $W(k)$-alg\`ebres locales noeth\'eriennes compl\`etes de corps
r\'esiduel $k$. Nous notons $\mathcal{M}_A$ l'id\'eal maximal d'un objet $A$ de $\widehat{\CC_k}$. Soit $X$ une courbe projective lisse sur $k$ et $G$ un groupe
fini d'automorphismes de $X$. Notons $D(G,X)$ le foncteur des d\'eformations
\'equivariantes de $(X,G)$ de $\CC_k$ dans la cat\'egorie des ensembles  qui
\`a un objet $R$ de $\CC_k$ associe l'ensemble des rel\`evements de $(X,G)$
\`a $R$ \`a isomorphisme \'equivariant pr\`es.
Nous savons
 que $D(X,G)$ admet un anneau versel $R(X,G)$ de d\'eformations
(\'equivariantes)  dans $\widehat{\CC}_k$. Rappelons
 que $\nu(X,G)=0$ signifie qu'il existe un rel\`evement
de $(X,G)$ en caract\'eristique z\'ero. Par ailleurs $\nu(X,G)=p^n$ 
 signifie qu'il existe un rel\`evement
de $(X,G)$ \`a un objet $A$ de $\CC_k$ de caract\'eristique $p^n$ et que pour
tout objet $A'$ de $\CC_k$ de caract\'eristique $p^{n+1}$,
$D(G,X)(A')=\emptyset$.\\

\noindent
D'apr\`es le principe local-global, l'\'etude de $D(X,G)$ se localise aux
points de ramification $\{P_1,\ldots,P_r\}\subset X$ :
$$R(X,G)\simeq R_1\hat\otimes\cdots\hat\otimes R_r[[U_1,\ldots,U_N]]$$
o\`u $R_i$ est l'anneau versel du foncteur $D_{G_{P_i}}$ des d\'eformations
infinit\'esimales de l'action du stabilisateur $G_{P_i}$ de $P_i$ sur le
compl\'et\'e de l'anneau local de $X$ en $P_i$ ($1\leq i\leq r$)
et $N$ est la dimension de l'espace des d\'eformations localement triviales
(voir Corollaire 3.3.5 \cite{BeMe}). Le calcul de la caract\'eristique
$\nu(X,G)$ se ram\`ene donc au calcul des caract\'eristiques des anneaux
$R_i$.\\

\noindent Dans la partie suivante, on va d\'emontrer le th\'eor\`eme ``local'' suivant (voir Propositions
\ref{pgen}, \ref{prop2} et \ref{fin}) :
\begin{theo}
\label{thlocal}
Soit $k$ un corps parfait de caract\'eristique $p>0$ et
$G\subset\Aut k[[y]]$ une action faible \textup{(}$G_2=1$\textup{)}.\\
L'action de $G$ se rel\`eve en caract\'eristique z\'ero si
$G$ est de l'une des trois formes suivantes :\\
\indent \textup{i.} groupe cyclique $\Z/p\Z$ d'ordre $p$ ;\\
\indent \textup{ii.} groupe di\'edral $D_p=\Z/p\Z\rtimes\Z/2\Z$ d'ordre $2p$ si $p>2$
ou $(\Z/2\Z)^2$ si $p=2$ ;\\
\indent \textup{iii.} groupe t\'etra\'edrique $A_4=(\Z/2\Z)^2\rtimes\Z/3\Z$ pour $p=2$.\\
Dans tous les autres cas, l'action de $G$ ne se rel\`eve pas en
caract\'eristique $p^2$; donc l'anneau $R_G$ versel de d\'eformations de l'action locale de $G$ satisfait $\car(R_G) \in \{0,p\}$ dans tous les cas. \blok
\end{theo}

D'apr\`es le principe local-global,  le th\'eor\`eme \ref{thlocal} donne
le r\'esultat global annonc\'e dans l'introduction (th\'eor\`eme \ref{theo1}).

\begin{rem}
Les d\'emonstrations impliquent {\sl a posteriori} que tout re\-l\`e\-ve\-ment d'une action faiblement ramifi\'e sur
$k[[y]]$ se rel\`eve (\`a conjugaison pr\`es) par homographies en $y$. Nous ne connaissons pas de d\'emonstration directe de cette
propri\'et\'e, sauf pour les rel\`evements \`a un anneau hens\'elien \'equicaract\'eristique (voir \cite{CoKa2}, 1.18 et Anhang B). 
\end{rem}

\section{\'Etude locale}

Dans la suite du paragraphe \S 3, nous nous pla\c{c}ons en un point $P\in X$
 de ramification sauvage (le cas d'une action mod\'er\'ee est bien connu,
 \cite{Gr}). Pour all\'eger les notations, \'ecrivons $G=G_P$. Nous 
supposons que
 l'action est faible : $G_{2}=1$. Il existe donc $t\geq 1$ et $n$ divisant $p^t-1$ tels que
$$G=(\Z/p\Z)^t\rtimes\Z/n\Z$$
avec une action de conducteur 1  (voir
\cite{CoKa2} \S 1, \cite{Se} VI.2).
Notons $k[[y]]$ l'anneau local compl\'et\'e de la courbe $X$ au point $P$
et 
$M$ la fibre compl\'et\'ee du faisceau tangent au point $P$.
L'espace tangent au foncteur de d\'eformations infinit\'esimales $D_G$ est
$H^1(G,M)$. La fibre $M$ avec action de $G$ s'identifie au $k[[y]]$-module des champs de vecteurs formels
$k[[y]]\frac{d}{dy}$ (\cite{BeMe} \S 2). 

\subsection{L'action \'equicaract\'eristique de $G$}

Nous commen\c{c}ons
 par supposer
que $n=1$. A conjugaison pr\`es par un automorphisme continu de la cl\^oture
alg\'ebrique $\overline{k((y))}$, l'action de $G$ est de la forme
\begin{equation}
\label{equa1}
G = \{ \sigma_u \}_{u \in V} \mbox{ avec } \sigma_u(y)= \frac{y}{1+u y}
\end{equation}
o\`u $V$ d\'ecrit un espace vectoriel $V \subseteq k$ de dimension $t$ sur $\F_p$ (voir \cite{CoKa} ou \cite{CoKa2}). Deux repr\'esentations conjugu\'ees de $G$ dans
$\Aut_kk[[y]]$ ayant m\^eme anneau versel de d\'eformations
(infinit\'esimales), nous nous ramenons au cas o\`u l'action de $G$ est de la
forme
(\ref{equa1}). Autrement dit l'action de $G$ est donn\'ee, \`a \'equivalence pr\`es, par un plongement
$G\hookrightarrow\PGL_2(k)$. \\

\noindent Un rel\`evement de $G$ \`a un anneau de $\CC_k$ de caract\'eristique
quelconque donne par r\'eduction modulo $p$
un rel\`evement \'equicaract\'eristique. Nous connaissons la d\'eformation verselle \'equicaract\'eristique
 (voir \S \ref{ck}). 
Il s'agit donc de calculer
explicitement les rel\`evements en caract\'eristique mixte de cette d\'eformation verselle. On utilisera la connaissance de l'anneau versel de d\'eformation d'un groupe cyclique (voir \S \ref{cycl}).  
Le cas $n$ quelconque (\S \ref{cas3})
d\'ecoule naturellement du cas $n=1$ (\S\S \ref{cas1} -- \ref{cas2}).
\begin{rem}
Dans cet article, les \'egalit\'es matricielles, traduisant
des \'egalit\'es entre homographies, doivent \^etre lues dans PGL$_2$, c'est-\`a-dire \`a multiplication
par un scalaire non nul pr\`es.
\end{rem}
\subsection{Anneau de d\'eformations \'equicaract\'eristiques}\label{ck} Rappelons les r\'esultats principaux
de \cite{CoKa} \S 4 : 
\begin{prop}
Soit $$G=(\Z/p\Z)^t \cong \{ \sigma_u = \left( \begin{array}{cc} 1 & 0 \\ u & 1 \end{array} \right) : u \in V
\}, 
 \mbox{ pour } {\displaystyle  V = \bigoplus_{i=1}^t {\bf F}_p \cdot u_i}.$$
Si $p\geq 3$ et $t>1$, l'anneau versel de d\'eformations \'equicaract\'eristiques de $G\subset\Aut k[[y]]$ est

$$ R_G/p = k[[\alpha,x_1,\ldots,x_{t-1}]]/\langle\alpha^{\frac{p-1}{2}},\alpha x_1,\ldots,\alpha
x_{t-1} \rangle.$$
La d\'eformation verselle \'equicaract\'eristique est la classe d'isomorphisme du rel\`evement d\'efini par $\{ \tilde{\sigma}_u \}_{u \in V}$ o\`u
$$
 {\displaystyle \tilde{\sigma}_u = \left(
\begin{array}{cc}
\displaystyle \sum^{\frac{p-1}{2}}_{j=0}{u+j-1\choose 2j}\alpha^j&
\displaystyle\alpha\sum^{\frac{p-1}{2}-1}_{j=0}{u+j\choose 2j+1}\alpha^j\\
\displaystyle\sum^{\frac{p-1}{2}-1}_{j=0}{u+j\choose 2j+1}\alpha^j+u_1+\sum_{i=1}^{t-1} x_i u_i&
\displaystyle\sum^{\frac{p-1}{2}}_{j=0}{u+j\choose 2j}\alpha^j
\end{array}
\right)},
$$
et pour $x\in k$ et $n\in\{0,\ldots,p-1\}$, le coefficient $\binom{x}{n}$ est d\'efini par $$\binom{x}{n}=\frac{x(x-1)\cdots(x-n+1)}{n!}.$$
\end{prop}
\begin{rem}
Si $p=3$, le probl\`eme est rigide ($R_G/3=k$) et $\alpha=0$.
Si $p=2$, le r\'esultat est plus compliqu\'e, mais ne sera pas utilis\'e ici. 
\end{rem}
\begin{rem}
Dans les d\'emonstrations, on utilisera surtout que les matrices $\tilde{\sigma}_u$ sont de la forme $\left( \begin{array}{cc} A & \alpha C \\ C & A+\alpha C \end{array} \right)$ pour certains $A,C \in R_G/p$.
\end{rem}
\subsection{Anneau de d\'eformations d'un groupe cyclique}\label{cycl} Rappelons les r\'esultats connus pour
$G=\Z/p\Z$ (\cite{BeMe}, th\'eor\`eme 4.2.8) :
\begin{prop}
\label{bm}
Soit $G=\Z/p\Z=\langle \sigma \rangle \subset\Aut k[[y]]$ avec $\displaystyle\sigma(y)=\frac{y}{y+1}$.\\
\indent \textup{i.} Si $p\geq 3$, l'anneau versel de d\'eformations infinit\'esimales de $G\subset\Aut k[[y]]$ est
$R_G=W(k)[[\alpha]]/\psi(\alpha)$ avec
\begin{equation} \label{psi} \psi(\alpha)=\sum_{\ell=0}^{\frac{p-1}{2}}\binom{p-1-\ell}{\ell}(-1)^{\ell}(\alpha+4)^{\frac{p-1}{2}-\ell}.\end{equation}
La d\'eformation verselle est la classe d'isomorphisme du rel\`evement d\'efini par
$$\displaystyle\widetilde{\sigma}(y)=\frac{y+\alpha}{y+\alpha+1}.$$En particulier, si $p=3$, $\psi(\alpha)=\alpha+3$, $R_G=W(k)$, le probl\`eme de d\'eformation est rigide et la d\'eformation verselle
est d\'efinie par
$$\displaystyle\widetilde{\sigma}(y)=\frac{y-3}{y-2}.$$

\indent \textup{ii.} Si $p=2$, $R_G=W(k)[[\alpha]]$ et la d\'eformation verselle
est d\'efinie par
$$\displaystyle\widetilde{\sigma}(y)=\frac{y+\alpha}{y-1}.$$
\end{prop}

\subsection{Cas $n=1,p \geq 3$}\label{cas1}
Soit $G=(\Z/p\Z)^t$ avec $t\geq 2$ et $p \geq 3$. 

\begin{lem}
\label{lgen}
Soit $H=(\Z/p\Z)^2=\langle \sigma,\tau \rangle \subset \Aut k[[y]]$
un $p$-groupe commutatif ($p \geq 3$), o\`u
$$\sigma(y)=\frac{y}{1+y},\;\; \tau(y)=\frac{y}{uy+1},$$
pour $u\in k-\F_p$. Alors l'anneau versel  $R_H$ des d\'eformations infinit\'esimales de $H\subset \Aut
k[[y]]$ est de caract\'eristique $p$.
\end{lem}
\begin{pr} D'apr\`es \S \ref{cycl}, $W(k)[[\alpha]]/(\psi(\alpha))$ est l'anneau versel de 
d\'eformation du sous-groupe $\langle\sigma\rangle\subset\Aut k[[y]]$ et la d\'eformation verselle de
$\langle\sigma\rangle\subset\Aut k[[y]]$ est donn\'ee par la classe d'isomorphisme du rel\`evement d\'efini par
$\widetilde{\sigma}(y)=\frac{y+\alpha}{y+\alpha+1}$. Il existe donc un morphisme 
$$W(k)[[\alpha]]/(\psi(\alpha))\rightarrow R_G$$
qui induit la d\'eformation du sous-groupe $\langle\sigma\rangle$ \`a $R_G$. Par cons\'equent la d\'eformation verselle
de $H\subset \Aut k[[y]]$ \`a $R_H$ est donn\'ee sur $\langle\sigma\rangle$ par le rel\`evement
$$m(y)=\frac{y+\alpha'}{y+\alpha'+1},$$
pour $\alpha'\in R_H$ avec $\psi(\alpha')=0$, et on denote encore $\alpha'$, par abus de notation, par $\alpha$.\\
De plus, d'apr\`es \S \ref{ck}, la r\'eduction modulo $p$ de la restriction \`a $\langle\tau\rangle$
de la d\'eformation verselle de $H$ \`a $R_H$ est d\'efinie par
$$\bar{n}(y)=\frac{Ay+\tilde{\alpha}C}{Cy+A+\tilde{\alpha}C},$$
o\`u $A$ et $C$ sont les polyn\^omes en $\tilde{\alpha}$ et $u$ d\'efinis au \S \ref{ck} 
et $\tilde{\alpha}\in R_H/p$.\\

\noindent {\sc \'Egalit\'e des param\`etres des deformations.} \ 
Nous commen\c{c}ons par d\'emontrer que $\tilde{\alpha}\equiv \alpha$ dans $R_H/p$. D'apr\`es \S \ref{ck}, la
d\'eformation verselle \'equicaract\'eristique est donn\'ee par des rel\`evements de la forme
$$\widetilde{\sigma}_v=\left(\begin{array}{ccc}A&\tilde{\alpha}C\cr C&A+\tilde{\alpha}C
\cr\end{array}\right).$$
Or l'image $\overline{m}$ de $m$ dans $R_H/p$ est de la forme
$$\left(\begin{array}{ccc}1&\alpha\cr 1&1+\alpha
\cr\end{array}\right).$$
Comme cette matrice doit \^etre \'egale \`a une des matrices $\tilde{\sigma}_v$ dans $PGL(2)$, on trouve que pour un tel $\tilde{\sigma}_v$, $A=C$ et alors aussi $\alpha \equiv \tilde{\alpha}$ dans $R_H/p$.\\

\noindent
Supposons $R_H\not=R_H/p$ et notons $R=R_H/p^2$. Notons $$n(y)=\frac{Ay+\alpha C}{Cy+ \alpha C+A}$$ avec $A,C$ les rel\`evements triviaux de $A$ et $C$ \`a $R$. Soit 
$$T(y):=n(y)+pS(y) \mbox{ pour } S \in R[[x]]$$
un rel\`evement arbitraire de $\bar{n}(y)$ \`a $R[[x]]$. \\
Il s'agit de voir que pour tout $S\in R[[x]]$, $\langle m,T \rangle \subseteq \Aut R[[x]]$ ne d\'efinit pas un rel\`evement du \emph{groupe} $H$. Pour parvenir \`a une contradiction, on va exprimer les relations dans $H$ en terme de la s\'erie $S$ tronqu\'ee \`a un certain ordre en $y$. \\

\noindent {\sc Relation de commutation et \'equation fonctionnelle.} \ La relation de commutation $\sigma\circ\tau=\tau\circ\sigma$ se traduit dans $R[[y]]$ par
$$m\circ(n+pS)=(n+pS)\circ m.$$
D'o\`u
$$\frac{n(y)+pS(y)+\alpha}{n(y)+pS(y)+\alpha+1}=n(m(y))+pS(m(y)).$$
Or un calcul direct montre que $n\circ m=m\circ n$ dans $R[[y]]$, soit encore
$$\frac{n(y)+pS(y)+\alpha}{n(y)+pS(y)+\alpha+1}=\frac{n(y)+\alpha}{n(y)+\alpha+1}+pS(m(y)).$$
D'o\`u l'\'equation fonctionnelle satisfaite par $S$
$$S(m(y))=(n(y)+\alpha+1)^{-2}S(y),$$
dans $R/p[[y]]$.\\

\noindent {\sc Solution tronqu\'ee de l'\'equation fonctionnelle.} \ On continue les calculs en utilisant une racine de $\psi$ dans une extension de $W(k)$ (projet\'ee dans $R$). Soit donc $K$ le corps des fractions de $W(k)$. Soit $L$ le corps de d\'ecomposition
de $\psi(X)$ sur $K$. Soit $\pi$ un id\'eal premier de $L$ au-dessus de
  $p$. Soit $L_{\pi}$ le localis\'e de $L$ en $\pi$, d'anneau des entiers
$\OO$. Soit $a$ une racine de $\psi$ dans $L_\pi$. Comme $\psi$ est unitaire,
$a\in\OO$. On observe que $p$ est ramifi\'e dans $L$, car $\psi(X) = X^{\frac{p-1}{2}}$ mod $p$. 
 L'\'equation fonctionnelle de $S$ mod $p$ implique la m\^eme \'equation pour $S$ mod $\pi$, 
  $\alpha=a=0$ mod $\pi$ ; donc
\begin{equation} \label{eqfonct2} S\Big(\frac{y}{y+1}\Big) = \Big(\frac{y}{uy+1}+1\Big)^{-2} 
\cdot S(y) \mbox{ dans } R/\pi[[y]]. \end{equation}
\noindent
La substitution de $S(y) = c_0 + c_1 y + c_2 y^2 + O(y^3)$ dans cette \'equation donne les formules suivantes
:
$$\left\{\begin{array}{cc} 2c_0=0,\cr
 c_1-(2u+3)c_0=0.\cr \end{array}\right.$$
Comme $p \neq 2$, nous obtenons $c_0=c_1=0$. 
Nous pouvons donc supposer que $S(y) = c y^2$ mod $(\pi,y^3)$ pour une constante $c$ convenable. \\

\noindent {\sc Lin\'earisation de la condition d'ordre $p$.} \ 
Nous continuons \`a \'ecrire $n(y)$ pour $n(y)|_{\alpha=a}$. Exprimons \`a pr\'esent
 que $T(y):=n(y)+pS(y)$ est d'ordre $p$ dans $R$. Montrons par r\'ecurrence sur $0\leq i\leq p$ que :
\begin{equation} \label{ind} T^i(y) \equiv n^i(y) + picy^2 \mbox{ mod } (p\pi,y^3). \end{equation}
C'est vrai pour $i=1$. \'Ecrivons 
$$n^k(y) = \frac{A_ky+B_k}{C_ky+D_k}.$$
Alors, modulo $y^3$ on trouve :
\begin{eqnarray*} T^{i+1}(y) &\equiv& T^{i}(T(y))\equiv n^{i}(T(y)) + p i c (n(y)+pc y^2)^2 \\ 
&\equiv& n^{i}(n(y)+pcy^2)+pic n(y)^2 \\ &\equiv& \frac{A_{i}n(y)+B_i+A_i p c y^2}{C_i n(y) + D_i+pcC_i 
y^2} + pic y^2 \mbox{ mod } (p\pi,y^3). \end{eqnarray*}
Utilisons l'identit\'e $$\frac{1}{r+ps} \equiv \frac{1}{r} - p \frac{s}{r^2} \mbox{ mod } p \pi $$ pour d\'evelopper
le d\'enominateur :
\begin{eqnarray*} T^{i+1}(y) &\equiv &  \frac{A_{i}n(y)+B_i}{C_i n(y) + D_i}+\frac{A_i p c y^2}{C_i n(y) + D_i} \\ & & -pc
C_i y^2 \frac{A_i n(y)+B_i}{(C_i n(y) + D_i)^2}+ pic y^2 \mbox{ mod } (p\pi,y^3). \end{eqnarray*} 
Modulo $p\pi$, toute expression de la forme $pX(a)$ peut \^etre remplac\'ee par $pX(0)$ : 
\begin{eqnarray*} T^{i+1}(y) &\equiv& n^{i+1}(y) + pc y^2 ( \frac{uy+1}{(i+1)uy+1} - 
\frac{iuy(uy+1)}{((i+1)uy+1)^2}+i) \\ &\equiv& n^{i+1}(y) + pc(i+1)y^2 \mbox{ mod } (p\pi,y^3), 
 \end{eqnarray*}
et (\ref{ind}) est d\'emontr\'ee. Ceci implique en particulier que 
$$ T^p(y) \equiv n^p(y) \mbox{ mod } (p \pi, y^3), $$ et la condition que $T^p(y)=y$ se traduit par 
la condition matricielle 
$$n^p(y) \equiv {\bf 1} \mbox{ mod } (p \pi, y^3).
$$ Observons maintenant que pour une homographie $H$ arbitraire $$ H(y):=\frac{ay+b}{cy+d} \mbox{ on a } H(y) = \frac{b}{d} + \frac{ad-bc}{d^2}
\cdot y + \frac{c(bc-ad)}{d^3} y^2 + O(y^3), $$
donc  $H(y)=y$ mod $y^3$ \'equivaut \`a $b=0, a=d, c=0$, c.-\`a-d.\ $H(y)=y$. 
La condition matricielle est donc \`equivalente \`a
$$n^p(y) \equiv {\bf 1} \mbox{ mod } p \pi. $$ 

\noindent {\sc Diagonalisation de la matrice d'ordre $p$.} \ 
Nous allons r\'einterpr\'eter cette condition en diagonalisant la matrice $n$ (sur $L$).
La matrice $n= \left(\begin{array}{cc}A&a C\cr C&
  A+a C\cr\end{array}\right)$ a pour polyn\^ome caract\'eristique
$$X^2-(2A+a C)X+A^2+a AC-a C^2$$
de discriminant $a C^2(a+4)$. Les valeurs propres sont donc
$$\lambda_u^{\pm}=A-C+Ce^{\pm i\theta}$$ avec la notation formelle
$$e^{\pm i\theta}=\cos\theta+i\sin\theta=Y\pm\sqrt{Y^2-1}\mbox{ pour }
Y=\cos\theta=1+a/2.$$
Comme $a\not=4$, $n$ est diagonalisable. La matrice $n$ est d'ordre $p$ si et seulement si
$(\lambda_u^{\pm})^p=1$.
Ces conditions se traduisent comme $$P(a)=Q(a)=0,$$ pour
\begin{eqnarray*} P(X)&=&\sum_{j=0}^p\binom{p}{j}(A-C)^{p-j}C^jT_j(X)-1 \\ Q(X)&=&\sum_{j=0}^p\binom{p}{j}(A-C)^{p-j}C^jS_{j-1}(X)\end{eqnarray*}
o\`u $T_j$ et $S_{j-1}$ sont les polyn\^omes de Tchebychef respectivement de
premi\`ere et de seconde esp\`ece d\'efinis par 
$e^{i j \theta} = T_j(\cos \theta) + i \sin \theta S_{j-1} (\cos \theta)$ et donn\'es explicitement par les formules :
\begin{eqnarray*} T_j(X)&=&\frac{1}{2}\sum_{\ell=0}^{\lceil
  j/2\rceil}\binom{j-\ell}{\ell}\frac{j}{j-\ell}(-1)^{\ell}(2X)^{j-2\ell}\\
S_{j-1}(X)&=&\sum_{\ell=0}^{\lceil
  (j-1)/2\rceil}\binom{j-1-\ell}{\ell}(-1)^{\ell}(2X)^{j-1-2\ell}\end{eqnarray*}

\noindent {\sc Calcul de la relation explicite d'ordre $p$.} \ Calculons $Q(a)$ modulo $p\pi$. Nous constatons que
$S_{-1}=0$, $S_{p-1}(a)=0$, car $\psi(X)$ est, par d\'efinition, le
  g\'en\'erateur de l'id\'eal $$ \langle T_p(1+X/2)-1,S_{p-1}(1+X/2) \rangle,$$ \cite{BeMe}
  lemme 4.2.6. En outre, $p$ divise $\binom{p}{j}$ pour $1\leq j\leq p-1$. Pour
  calculer $Q(a)$ modulo $p\pi$, nous pouvons
 donc remplacer $A,C,(S_{j-1})_{0\leq
  j\leq p}$ par leurs valeurs modulo $\pi$. Or $S_{j-1}\equiv j\mod\pi$. Alors
 \begin{eqnarray*}
Q(a) &=& p (\sum_{j=1}^{p-1}\frac{j}{p}\binom{p}{j}(A-C)^{p-j}C^j) \mod p \pi \\ &=& p (\sum_{j=1}^{p-1} \binom{p-1}{j-1}(A-C)^{p-1-(j-1)}C^j) \mod p \pi \\
&=& pC (\sum_{k=0}^{p-2} \binom{p-1}{k}(A-C)^{p-1-k}C^k) \mod p \pi \\
&=& pC (A^{p-1}-C^{p-1}) \mod p \pi
\end{eqnarray*}
Si $Q(a)=0 \mod p \pi$, il faut que $C(A^{p-1}-C^{p-1})=0 \mod \pi$, d'o\`u $C=\zeta A$ pour un $\zeta \in \F_p$. Mais alors la matrice de $n$ est de la forme
$$ \left( \begin{array}{cc} A & a C \\ C & A+aC \end{array} \right) = A \left( \begin{array}{cc} 1 & a \zeta \\ \zeta & 1+a\zeta \end{array} \right) \mod p, $$
ce qui implique que  $C\mbox{ mod } p  \in {\bf F}_p$. Or $C \mbox{ mod } (p,\mathfrak{m}_R) \equiv u \notin \F_p$, 
ce qui conduit \`a une contradiction. \end{pr}

\begin{prop}
\label{pgen}
Soit $p \geq 3$. Supposons que l'action de $G=(\Z/p\Z)^t$, $t\geq 2$ sur $k[[y]]$ soit faiblement ramifi\'ee. Alors
l'anneau versel $R_G$ de d\'eformations infinit\'esimales  est de caract\'eristique $p$.
\end{prop}
\begin{pr}
Nous pouvons normaliser l'action de $G\subset\Aut k[[y]]$ de fa\c{c}on \`a ce qu'il existe $\sigma, \tau$ deux
\'el\'ements de $G$ tels que
$$\sigma(y)=\frac{y}{1+y},\;\; \tau(y)=\frac{y}{uy+1},$$
pour $u\in k-\F_p$. Soit $H$ le sous-groupe de $G$ engendr\'e par $\sigma$ et $\tau$. D'apr\`es le lemme
\ref{lgen}, l'anneau versel de d\'eformations infinit\'esimales $R_H$ est de caract\'eristique $p$.
Par versalit\'e de $R_H$, il existe un morphisme $R_H\rightarrow R_G$. Donc $R_G$ est de
caract\'eristique $p$. 
\end{pr}

\subsection{Cas $n=1,p =2$}\label{cas2}

\begin{prop}
\label{prop2}
Supposons $p=2$. Soit $t\geq 1$ et
$G=(\Z/2\Z)^t\subset \Aut k[[y]]$ de conducteur 1.\\
\indent \textup{i.} Si $t\leq 2$, alors l'action de $G$ se rel\`eve en caract\'eristique z\'ero, et $\car (R_G)=0$.\\
\indent \textup{ii.} Si $t\geq 3$, alors l'action de $G$ ne se rel\`eve pas en carat\'eristique
$p^2$, et $\car (R_G)=p$. 
\end{prop}
\begin{pr}
\indent \textup{i.} Par conjugaison, nous pouvons supposer que $H=\langle \sigma_1 \rangle$ est un sous-groupe de $G$. 
D'apr\`es \S\ref{cycl}, la d\'eformation verselle de $H\subset\Aut k[[y]]$ \`a $R_H=W(k)[[\alpha]]$ est
donn\'ee par
$$\tilde{\sigma}_1=\left(\begin{array}{cc}1&\alpha\cr 1&-1\cr\end{array}\right).$$ 
(D'o\`u i.\ pour $t=1$). Supposons $t=2$. Par conjugaison, nous pouvons supposer
qu'il existe $\sigma_u\in G$ d'ordre 2 tel que
$$\sigma_u(y)=\frac{y}{u+y}$$
avec $u\in k-\F_2$. Soit
$$\tilde{\sigma}_u=\left(\begin{array}{cc}1&-\alpha
  \tilde{u}-2\cr \tilde{u}&-1\cr\end{array}\right)$$
   pour $\tilde{u}$ un
  rel\`evement de $u$ \`a $W(k)$. Alors un tel $\tilde{\sigma}_u$ est la seule matrice qui satisfait $\tilde{\sigma}_u^2=\tilde{\sigma}_1^2=1$
et $\tilde{\sigma}_u\tilde{\sigma}_1=\tilde{\sigma}_1\tilde{\sigma_u}$. Par cons\'equent
 $G=(\Z/2\Z)^2$ se rel\`eve en
caract\'eristique 0 \`a $W(k)[[\alpha]]$.
Ce rel\`evement est versel car dim $H^1(G,M)=1$. \\

ii. Supposons $t\geq 3$. Soient $\sigma_1,\sigma_u,\sigma_v$ des \'el\'ements de $G\subset\Aut k[[y]]$
avec 
$$\sigma_1(y)=\frac{y}{1+y},\;\;\sigma_u(y)=\frac{y}{u+y},\;\;\sigma_v(y)=\frac{y}{v+y},$$
avec $\{1,u,v\}$ $\F_2$-lin\'eairement ind\'ependants.
Soit $R_G$ l'anneau versel de d\'eformations infinit\'esimales de $G$. Par versalit\'e, nous avons des
morphismes 
$$W(k)[[\alpha]]\rightarrow R_G,\;\; W(k)[[\beta]]\rightarrow R_G$$
qui induisent la d\'eformation verselle  \`a $R_G$ donn\'ee par les rel\`evements
$$\tilde{\sigma}_1=\left(\begin{array}{cc}1&\alpha\cr 1&-1\cr\end{array}\right)=
\left(\begin{array}{cc}1&\beta
  \cr 1&-1\cr\end{array}\right)$$
  $$\tilde{\sigma}_u=\left(\begin{array}{cc}1&-\alpha
  \tilde{u}-2\cr \tilde{u}&-1\cr\end{array}\right),\;\;\tilde{\sigma}_v=\left(\begin{array}{cc}1&-\beta
  \tilde{v}-2\cr \tilde{v}&-1\cr\end{array}\right)$$
avec $\alpha,\beta\in R_G$ (les \'egalit\'es sont dans $PGL(2)$). D'o\`u $\alpha=\beta$.
De plus les \'equations de commutation
$$\tilde{\sigma}_1 \tilde{\sigma}_u = \lambda_{u,1} \tilde{\sigma}_u \tilde{\sigma}_1, \ \tilde{\sigma}_1 \tilde{\sigma}_v =\lambda_{v,1} \tilde{\sigma}_v \tilde{\sigma}_1, \ \tilde{\sigma}_u \tilde{\sigma}_v = \lambda_{u,v} \tilde{\sigma}_v \tilde{\sigma}_u $$
sont satisfait exactement pour $\lambda_{u,1}=\lambda_{v,1}=\lambda_{u,v}=-1$ et $$2\alpha = -
2\frac{\tilde{u}+\tilde{v}-1}{\tilde{u}\tilde{v}} $$ mais en caract\'eristique 4, pour cette valeur de $\alpha$, $\tilde{\sigma}_v = \tilde{\sigma}_1 \tilde{\sigma}_u$. 
Il n'existe donc pas de rel\`evement de $(\Z/2\Z)^3$ en caract\'eristique 4.
\end{pr}
\subsection{Cas $n>1$}\label{cas3}
\begin{prop}
\label{fin}
Soit $G=(\Z/p\Z)^t\rtimes \Z/n\Z\subset\Aut k[[y]]$ une action faiblement ramifi\'ee \textup{(}$G_2=1$\textup{)} avec
$n>1$ et $n$ divisant $p^t-1$.\\
\indent \textup{i.} Si $t>1$, l'action de $G$ se rel\`eve en caract\'eristique $p^2$ si et
seulement si $p=2$ et
$$G=(\Z/2\Z)^2\rtimes\Z/3\Z.$$
\indent \textup{ii.} Si $t=1$,  l'action de $G$ se rel\`eve en caract\'eristique $p^2$ si et
seulement si $p\geq 3$ et
$$G=\Z/p\Z\rtimes\Z/2\Z. $$
De plus dans les cas o\`u l'action de $G$ se rel\`eve en caract\'eristique
$p^2$, elle se rel\`eve en caract\'eristique z\'ero, et alors $\car(R_G)=0$.
\end{prop}
\begin{pr}
\textup{i.} Supposons $t\geq 2$.\\
\noindent
 Si $p\not=2$, d'apr\` es les propositions \ref{pgen}, 
 le $p$-Sylow de $G$ ne se rel\`eve pas en 
caract\'eristique $p^2$.\\
Si $p=2$, d'apr\`es la proposition \ref{prop2}, on peut supposer $t=2$. Comme
$n>1$ divise $p^t-1$, on a $n=3$ et $G=(\Z/2\Z)^2\rtimes\Z/3\Z$. D'apr\`es \cite{CoKa}, le probl\`eme de d\'eformations 
\'equicaract\'eristiques de $G$ est rigide. Le 
rel\`evement versel de $\langle\sigma_1\rangle$
\`a $W(k)[[\alpha]]$ est donn\'e par la matrice
$$m=\left(\begin{array}{cc}1&\alpha\cr 1&-1\cr\end{array}\right).$$
Par des
  calculs matricielles faciles, nous trouvons un rel\`evement unique (donc versel) de l'action de $G$ \`a
$W(k)[j]/\langle j^2+j+1 \rangle$ donn\'e par les matrices
$$m= \left(\begin{array}{cc}1&3j-1\cr 1&-1\cr\end{array}\right),\;\;
m'=\left(\begin{array}{cc}1&4j-1\cr j&-1\cr\end{array}\right),\;\;
  g=\left(\begin{array}{cc}1&-2j-1\cr 0&1\cr\end{array}\right)$$
satisfaisant $m^2={m'}^2=g^3=1$, $mm'+m'm=1$, $gmg^{-1}=jm'$.
 D'o\`u (i).\\

\indent \textup{ii.} Supposons $G=\Z/p\Z\rtimes\Z/n\Z$ avec $n>1$ divisant $p-1$
(donc $p>2$).\\

\noindent
$\bullet$
 Si $n\not=2$ et $p>3$, la d\'eformation verselle \'equicaract\'eristique de l'action de $\Z/p\Z \rtimes \Z/n\Z$ est rigide et de la forme (\`a conjugaison pr\`es)
$$ \langle \left( \begin{array}{cc} 1 & 0 \\ 1 & 1 \end{array} \right)\rangle \rtimes \langle \left( \begin{array}{cc} \zeta & 0 \\ 0 & 1 \end{array} \right) \rangle , $$ avec $\zeta$ $n$-i\`eme racine d'unit\'e ( voir \cite{CoKa}, 4.4.5(i); la pr\'esence de $\zeta$ ``tue'' donc les d\'eformations de $\Z/p\Z$ dans les directions ``$\alpha$''). Soit 
$$ T(y) = \frac{y}{y+1} + p \cdot S(y) \mbox{ pour } S \in R[[y]] $$
un rel\`evement \`a $R:=R_G/p^2$ du g\'en\'erateur de $\Z/p\Z$ donn\'e. Supposons car$(R)=p^2$. 
On pose $S(y) = c_0 + c_1 y + c_2 y^2 + O(y^3)$ et on trouve par induction que pour $N$ arbitraire, 
\begin{eqnarray*} T^N(y) &=& Npc_0 + (1+Npc_1 - N(N-1)c_0) x \\ & &  + (Npc_2 - N - \frac{3}{2} N(N-1) pc_1 + \frac{1}{3} N(4N^2-9N+5)pc_0)x^2 \\ & & + O(x^3) \mbox{ mod } p^2. \end{eqnarray*}
Pour $N=p$ et $p \neq 3$, ceci implique $T^p(x)=x-px^2 \mbox{ mod } (p^2,x^3), $
qui n'est jamais congru \`a $ x$ dans $R$. Un tel rel\`evement $T$ d'ordre $p$ n'existe donc pas.

\medskip

\noindent
$\bullet$ Si $p=3$, alors $n=2$.
Le groupe $G=\Z/3\Z\rtimes\Z/2\Z$ se rel\`eve (de fa\c{c}on rigide
d'apr\`es \cite{BeMe} 4.2.8)  \`a $W(k)$ par les matrices
$$m= \left(\begin{array}{cc}1&-3\cr 1&-2\cr\end{array}\right),\;\;
g=\left(\begin{array}{cc}1&0\cr 1&-1\cr\end{array}\right)$$
telles que $m^3=g^2=gmgm=1$.\\

\medskip
\noindent
$\bullet$ Le groupe $G=\Z /p\Z\rtimes\Z/2\Z$ avec $p>3$ admet un rel\`evement versel
$$m= \left(\begin{array}{cc}1&\alpha\cr 1&1+\alpha\cr\end{array}\right),\;\;
g=\left(\begin{array}{cc}1&\alpha\cr 0&-1\cr\end{array}\right)$$
\`a $W(k)[[\alpha]]/\langle\psi(\alpha)\rangle$ pour lequel $m^p=g^2=gmgm=1$. 
Ce rel\`evement est versel, car dim $H^1(G,M)=1$.\end{pr}

\noindent Les propositions \ref{pgen}, \ref{prop2} et \ref{fin} d\'emontrent le th\'eor\`eme ``local'' \ref{thlocal} et donc aussi le th\'eor\`eme ``global'' \ref{theo1}.

\section{Corollaires, remarques et questions ouvertes}
La d\'emonstration du th\'eor\`eme \ref{thlocal} et les r\'esultats de 
\cite{Gr}, \cite{BeMe} et \cite{CoKa} donnent des informations plus pr\'ecises
sur les anneaux versels de d\'eformations infinit\'esimales pour les actions
faibles.
\begin{cor}
\label{coroloc}
Soit $k$ un corps parfait de caract\'eristique $p>0$. Soit
$G=(\Z/p\Z)^t\rtimes\Z/n\Z\subset\Aut k[[y]]$ une action faible ($G_2=1$)
avec $t\geq 0$ et $n$ divise
$p^t-1$. Notons $R_G$ l'anneau versel de d\'eformations infinit\'esimales de
$G$. Soient $\zeta$ une racine primitive $n$-i\`eme de l'unit\'e dans $k$,
$s=[\F_p(\zeta):\F_p]$, et $\psi(\alpha)\in W(k)[\alpha]$ d\'efini par 
$$\psi(\alpha)=\sum_{\ell=0}^{\frac{p-1}{2}}\binom{p-1-\ell}{\ell}(-1)^{\ell}(\alpha+4)^{\frac{p-1}{2}-\ell}.$$
Alors nous avons  la table suivante :
$$\begin{tabular}{|l|l|}
\hline
$G$&$R_G$\\
\hline
$\Z/n\Z,(n,p)=1$  &$W(k)$ \\
\hline
$\Z/p\Z, p\not=2,3$ & $W(k)[[\alpha]]/\langle\psi(\alpha)\rangle$\\
\hline
$\Z/3\Z, p=3$& $W(k)$\\
\hline
$(\Z/2\Z)^t,p=2, t\leq 2$& $W(k)[[\alpha]]$\\
\hline
$(\Z/p\Z)^t,  t>1, p\neq 2,3$ & $k[[\alpha,x_1,\ldots,x_{t-1}]]/\langle\alpha^{\frac{p-1}{2}},\alpha x_1,\ldots,\alpha
x_{t-1}\rangle$\\
\hline
$(\Z/3\Z)^t,p=3,t>1$& $k[[x_1,\ldots,x_{t-1}]]$\\
\hline
$\Z/p\Z\rtimes\Z/2\Z, p\not=2,3$ & $W(k)[[\alpha]]/\langle\psi(\alpha)\rangle$\\
\hline
$\Z/3\Z\rtimes\Z/2\Z$& $W(k)$\\
\hline
$(\Z/2\Z)^2\rtimes\Z/3\Z$& $W(k)[j]/\langle j^2+j+1\rangle$\\
\hline
$(\Z/p\Z)^t\rtimes\Z/n\Z,$&\\
$ \ \  t\geq 2, n\not=2 \mbox{ ou } p=2,3$ & $
k[[x_1,\ldots,x_{t/s-1}]]$\\
\hline
$(\Z/p\Z)^t\rtimes\Z/n\Z,$&
\\
$ \ \  t\geq 2 \mbox{ et } p\not=2,3$ & $k[[\alpha,x_1,\ldots,x_{t/s-1}]]/\langle \alpha^{\frac{p-1}{2}},\alpha x_1,\ldots,\alpha x_{t-1}
\rangle$\\
\hline
\end{tabular}$$
Enfin si $G=(\Z/2\Z)^t, t>2$ alors 
$$R_G=k[[\alpha,x_1,\ldots,x_t]]/\langle x_1+\cdots+x_t,u_1x_1+\cdots+u_tx_t,\alpha(x_iu_j-x_ju_i)_{1\leq
  i,j\leq t}\rangle$$ o\`u  $\{u_1,\ldots,u_t\}$ est une base  de $V$.
\end{cor}
\begin{pr}
Les anneaux universels de d\'eformations de caract\'eristique z\'ero ont \'et\'e calcul\'es dans les preuves des
propositions \ref{prop2} et \ref{fin}. Dans les autres cas, les anneaux universels de d\'eformations ont pour
caract\'eristique $p$ et ont \'et\'e calcul\'es par Cornelissen et Kato (\cite{CoKa} \S 4.4).
\end{pr}
\begin{example}
Pour $p \neq 2,3$, la d\'eformation verselle \`a $R_G$ de $G=(\Z/p\Z)^t$ est donn\'ee dans \ref{ck}. \end{example}
\begin{rem}
Oort, Sekiguchi et Suwa (\cite{OoSeSu}) 
ont d\'emontr\'e que l'action  de $G=\Z/p\Z$ se
rel\`eve
en caract\'eristique z\'ero (pour tout conducteur).
Pagot a trait\'e le cas $G=(\Z/2\Z)^2$ (\cite{Pa}).
Bouw et Wewers ont trait\'e le cas du groupe di\'edral
$D_p=\Z/p\Z\rtimes\Z/2\Z$ (\cite{BoWe}).
\end{rem}
 
\begin{rem}
Le th\'eor\`eme \ref{theo1} pr\'ecise notamment quelles sont les actions sur les courbes
ordinaires qui se rel\`event en caract\'eristique z\'ero. Il r\'epond donc \`a
une version plus forte de la conjecture d'Oort (\cite{Oo})
dans le cadre des rev\^etements
faiblement ramifi\'es.
\end{rem}

\begin{rem}
Bertin \cite{Be}, Green \cite{Gre} Green et Matignon \cite{GrMa} ont donn\'e des crit\`eres
n\'ecessaires \`a l'existence d'un rel\`evement \emph{en caract\'eristique z\'ero}
d'un groupe de la forme $(\Z/p\Z)^t$ de conducteur arbitraire; l'impossibilit\'e de relever est impliqu\'e par certaines congruences satisfaites par les  conducteurs de Hasse locaux: cette th\'eorie permet une d\'emonstration assez simple du fait que $\nu(X,G) \neq 0$ dans le cas faiblement ramifi\'e o\`u les groupes de ramification non-trivials d'ordre divisible par $p$ ne sont pas sur la liste du corollaire \ref{theo1}. Le th\'eor\`eme
\ref{theo1} impose que le conducteur vaut 1 mais discute l'existence de
rel\`evements \emph{en caract\'eristique $p^n$ arbitraire} interm\'ediaire.
\end{rem}

\begin{rem}
Le th\'eor\`eme \ref{theo1} est compatible au resultat de Chinburg, Harbater et
Guralnick (\cite{ChHaGu}): ils donnent une liste de groupes abstraits $G$ 
ayant un $p$-Sylow normal $P$ avec $G/P$ cyclique, admettant
 une action $G \hookrightarrow \mbox{Aut}_k(k[[x]])$ dont
 l'obstruction locale au rel\`evement de Bertin (\cite{Be}) est non nulle. 
Les groupes i.--iii.\ du corollaire \ref{theo1} sont bien 
dans le compl\'ementaire de cette liste.
\end{rem}

\begin{rem}
 Il serait int\'eressant de donner des exemples de courbes $X$ munies d'une
 action finie $G$ telles que la caract\'eristique $\nu(X,G)$ de l'anneau
 versel de d\'eformations \'equivariantes soit de caract\'eristique $p^n$ pour
 $1<n<\infty$. Pour une discussion de cette question, voir \cite{Co}. 
\end{rem}

\begin{rem}
Les d\'emonstrations dans cet article sont une version explicite d'une strat\'egie g\'en\'erale pour
\'etudier le foncteur des d\'eformations \'equivariantes de $(X,G)$. L'id\'ee
est de d\'evisser le groupe $G$ et de contr\^oler le comportement du foncteur
$D_G$ des d\'eformations par restriction \`a un sous-groupe et par passage au
quotient. Le passage \`a un sous-groupe $H$ est \'el\'ementaire car nous avons un
morphisme de foncteurs canonique
$D_G\rightarrow D_H$
et les morphismes induits au niveau des espaces tangents et des groupes
d'obstruction co\"\i ncident canoniquement aux applications de restriction
cohomologique. Le passage \`a un quotient est beaucoup plus subtile. Il est
d\'ecrit localement par la suite spectrale d'Hochschild-Serre.
Dans cet article, nous avons men\'e explicitement les calculs d'obstruction
aux recollements
de $\Z/p\Z$ \`a $\Z/p\Z$. Il serait int\'eressant de
mener ces calculs d'obstructions aux 
recollements en termes cohomologiques (voir \cite{BeMe2}).
\end{rem}

\end{document}